\documentclass{article}
\usepackage{amsmath,amssymb,theorem,xypic}
\xyoption{curve}

\theoremstyle{change}
\theorembodyfont{\itshape}
\newtheorem{thm} {Theorem}
\newtheorem{prop}{Proposition}
\newtheorem{lemma}{Lemma}
\newtheorem{cor}{Corollary}

\numberwithin{equation}{subsection}
\numberwithin{thm}{subsection}
\numberwithin{prop}{subsection}
\numberwithin{lemma}{subsection}
\numberwithin{defn}{subsection}
\numberwithin{cor}{subsection}

\newcommand{\demobox}{\vrule height6pt width6pt depth0pt}

\makeatletter
\renewcommand{\subsection}{\@startsection%
{subsection}{2}{0mm}{\baselineskip}{-1em}%
{\normalfont\normalsize\bfseries}}
\renewcommand{\subsubsection}{\@startsection%
{subsubsection}{3}{0mm}{\baselineskip}{-1em}%
{\normalfont\normalsize\textit}}
\makeatother

\newenvironment{demo}{\noindent{\it Proof.}}
{{\unskip\nobreak\hfil\qquad\demobox}
\parfillskip=0pt\par\medskip}

\newcommand{\sset}{\subseteq}
\newcommand{\tens}{\otimes}
\newcommand{\Xto}{\xrightarrow}
\newcommand{\isom}{\Xto{\sim}}
\newcommand{\up}[2]{{}^{#1}{#2}}
\newcommand{\Ver}{\mathrm{Ver}}

\newcommand{\Sh}{\mathrm{Sh}}
\newcommand{\ad}{{\mathrm{ad}}}

\newcommand{\arith}{\mathit{arith}}
\newcommand{\ab}{\mathit{ab}}
\newcommand{\nr}{\mathit{nr}}
\renewcommand{\O}{\mathcal{O}}
\newcommand{\bQ}{\mathbb{Q}}
\newcommand{\bZ}{\mathbb{Z}}
\newcommand{\ctens}{\mathop{\hat\otimes}}
\newcommand{\Br}[1]{{{\operatorname{Br}}({#1})}}
\newcommand{\End}{{\operatorname{End}}}
\newcommand{\Ext}{{\operatorname{Ext}}}
\newcommand{\Gal}{{\operatorname{Gal}}}
\newcommand{\Hom}{{\operatorname{Hom}}}

\title{Weil groups and $F$-isocrystals}
\author{Richard Crew}
\date{}

\begin{document}
\maketitle

\section*{Introduction}
\label{intro}

The subject of this paper is an explicit formula for the fundamental
class $u_{L/K}$ of a finite Galois extension $L/K$ of local fields
with group $G$. This formula appears implicitly in in \cite[Ch. XIII
\S5 Ex. 2]{serre:1968} (see the discussion just before theorem
\ref{thm:dwork-serre-formula}), where Serre sketches a proof that this
class is indeed the fundamental class; his argument uses the theory of
class formations (and thus depends on the basic results of local class
field theory). In this paper we give a proof of this assertion
independent of any results of class field theory; we use instead the
Dieudonn\'e-Manin structure theorem for $F$-isocrystals over an
algebraically closed field. The formula yields a direct proof of a
generalization of Dwork's formula for the norm residue symbol (this is
also in \cite{serre:1968}) and this in turn yields a simple proof that
the cup product with $u_{L/K}$ is an isomorphism. In other words the
cohomological treatment of local class field theory can dispense with
the Tate-Nakayama theorem. 

As a by-product of our main construction we get a construction of the
relative Weil group $L/K$ as an automorphism group, and a simple proof
of the local Weil-Shafarevich theorem. When $K=\bQ_p$ this
construction of $W_{L/K}$ is due to Morava \cite{morava:1979}, who
expressed it in terms of Lubin-Tate groups. This construction leads
directly to a direct proof of Shafarevich's theorem (i.e.\ the local
Shafarevich-Weil theorem), which answers a question raised in
\cite{morava:1979}.

Our sign conventions for homological algebra are those of Bourbaki
\cite{bourbaki:hom}. We will follow Deligne \cite{deligne:1972} and
Tate \cite{tate:arcata} in normalizing the reciprocity law so that in
an unramified extension the arithmetic Frobenius corresponds to the
\textit{inverse} of a uniformizer; this makes the formulas come out
slightly simpler. We use Serre's normalization \cite[\S1
App.]{serre:1967} of the invariant of a division algebra over a local
field, which is the opposite of that used in \cite{demazure:1972}.

\textit{Acknowledgements.} I would like to thank the referee for a
thorough critical reading of the manuscript which resulted in many
improvements. A number of revisions were made during a visit to the
Kavli IPMU in the fall of 2024, which I would like to thank for its
hospitality.

\section{The Fundamental Class}
\label{sec:fundamental-class}

We first recall a few general results about $F$-isocrystals on a
field. In the next section we use these to construct the fundamental
class of of a finite Galois extension of local fields.

\subsection{$F$-isocrystals.}
\label{sec:F-isocrystals}

In this section (and \emph{only} in this section) $K$ will be a
complete nonarchimedean discretely valued field with algebraically
closed residue field $k$ of characteristic $p>0$. We do \textit{not}
assume that $K$ itself has characteristic zero. We fix a power $q$ of
$p$ and assume that $K$ has a lifting $\sigma$ of the $q$th power
Frobenius of $k$. Since $K$ is complete, $\sigma$ fixes some
uniformizer $\pi$ of $K$, which we also fix. Let $K_0$ denote the
fixed field of $\sigma$, so that $\pi\in K_0$.

An \textit{$F$-isocrystal on $K$} is a $K$-vector space $V$ of finite
dimension endowed with a $\sigma$-linear automorphism $F$, called the
\textit{Frobenius structure} of $V$ (or of $(V,F)$). This is the usual
notion when $K$ has characteristic $0$ but we will use this term when
$K$ has characteristic $p$ as well. If it is necessary to specify
$\sigma$ (and $q$) we will say $(\sigma,F)$-isocrystal. A morphism
$(V,F)\to(V',F')$ of $F$-isocrystals is a $K$-linear map $V\to V'$
compatible with $F$ and $F'$. Some authors refer to $F$-isocrystals
simple as isocrystals, but we will not do this.

Let $K_\sigma[F]$ be the Dieudonn\'e ring, i.e. the noncommutative
polynomial ring in $F$ with coefficients in $K$ and multiplication
defined by $Fa=a^\sigma F$. Evidently the category of $F$-isocrystals
on $K$ is equivalent to the subcategory of left $K_\sigma[F]$-modules
whose objects are $K$-vector spaces of finite dimension on which $F$
acts invertibly, and whose morphisms are $K_0$-linear maps compatible
with the actions of the $F$. 

The following theore gives the structure of this category:

\begin{thm}[Dieudonn\'e, Manin]
  The category of $F$-isocrystals on $K$ is semisimple, with one
  isomorphism class of simple objects for every element of $\bQ$. The
  class of $\lambda\in\bQ$ is the class of 
  \begin{equation}\label{eq:simple-F-isocrystal} 
    V_K(\lambda)=K_\sigma[F]/K_\sigma[F](F^d-\pi^r)
  \end{equation}
  where $\lambda=r/d\in\bQ$ in lowest terms (when $\lambda=0$ we take
  $r=0$ and $d=1$).
\end{thm}

In the original references this was proven under the assumption that
$K$ is of characteristic zero and absolutely unramified; see for
example Manin \cite[Ch. II]{manin:1963} or Demazure \cite[Ch. IV
\S4]{demazure:1972}. The only change needed for the general case is to
replace the prime $p$ by $\pi$ throughout; the original arguments work
since $\pi$ is fixed by $\sigma$.

The rational number $\lambda$ in \ref{eq:simple-F-isocrystal} is the
\emph{slope} of the $F$-isocrystal $V_K(\lambda)$. More generally, the
\emph{multiplicity of $\lambda\in\bQ$} as a slope of an $F$-isocrystal
$M$ is the $K$-dimension of the sum of the sub-$F$-isocrystals of $M$
isomorphic to $V_K(\lambda)$, and the \emph{slopes} of $M$ are the the
slopes of the $V_K(\lambda)$ occurring in $M$, each counted with its
multiplicity. Thus $V_K(\lambda)$ itself is purely of slope $\lambda$,
with multiplicity $d$.

We denote by $D_K(\lambda)$ the endomorphism algebra of
$V_K(\lambda)$. 

\begin{thm}[Dieudonn\'e]\label{thm:end-ring-of-simple-Fisocrystal}
  The $K_0$-algebra $D_K(\lambda)$ is a simple division algebra over
  $K_0$ with invariant $-\lambda$.
\end{thm}

When $K$ has characteristic $0$ and is absolutely unramified a proof
can be found in \cite[Ch. IV \S3]{demazure:1972}; in the general case
the changes mentioned above suffice. Note that \cite{demazure:1972}
has $\lambda$ for the invariant, since Demazure is using the opposite
convention for normalizing the invariant.

Thus if $\lambda=r/d$ in lowest terms, $D_K(\lambda)$ has dimension
$d^2$ over $K_0$. In fact it is not hard to check that $D_K(\lambda)$
is a cyclic algebra over $K_0$, but Dieudonn\'e's original proof works
just as well.

To compute the slopes of a simple $F$-isocrystal we will use the
following proposition, which is an easy consequence of Katz's ``basic
slope estimate'' \cite{katz:1979}. Let $\O_K$ be the integer ring
of $K$.  As always, a \textit{lattice} in a (finite-dimensional)
$K$-vector space $V$ is a finite $\O_K$-module $V_0$ such that
$K\otimes_{\O_K}V_0\simeq V$.

\begin{prop}\label{prop:bounds-on-Frobenius}
  The slopes of an $F$-isocrystal $(V,F)$ are contained in the
  interval $[\lambda,\mu]$ if and only if for some lattice
  $V_0\subset V$ there are constants $C$, $D$ such that
  \begin{equation}
    \label{eq:bounds-on-Frobenius}
    \pi^{[n\mu]+C}V_0\subseteq F^n(V_0)\subseteq\pi^{[n\lambda]+D}V_0
  \end{equation}
  independently of $n$.
\end{prop}

One can also prove this by observing first that if
\ref{eq:bounds-on-Frobenius} holds for one lattice it holds for any
lattice; this reduces to the case where $V$ is one of the simple
$F$-isocrystals \ref{eq:simple-F-isocrystal}, where it can be proven
by direct computation.

Note that for any $n>0$ the slopes of $(V,F)$ are determined by the
slopes of the $(\sigma^n,F)$-isocrystal $(V,F^n)$. We will use this
trick in the section \ref{sec:lubin-tate-fisoc}.

\subsection{The algebra $L_K$.}
\label{sec:L_K}

From now on $K$ will be a non-archimedean \emph{local field} with
residue field $k$ of characteristic $p$ and cardinality $q$. We denote
by $\O_K$ the integer ring of $K$, and use similar notation for all
non-archimedean fields.  Fix a completion $K^\nr$ of the maximal
unramified extension of $K$; its residue field is an algebraic closure
$\bar k$ of $k$. We write $v:K^\times\to\bZ$ for the normalized
valuation of $K$, i.e.\ $v(\pi)=1$ for any uniformizer $\pi$ of
$\O_K$. We will also use $v$ to denote the unique extension of $v$ to
any finite extension of $K$, to $K^\nr$ and finally to any finite
extension of $K^\nr$.

We start with some basic facts; the following two lemmas are probably
well-known:

\begin{lemma}\label{lemma:totally-ramified-extensions}
  Suppose $K$ is a complete non-archimedean field, $L/K$ is a finite
  totally ramified extension and $E/K$ is a discretely valued
  unramified extension of $K$. Then $L\tens_KE$ is a discretely valued
  field and the canonical homomorphism $E\to L\tens_KE$ induces an
  isomorphism of residue fields. If $E$ is complete, so is
  $L\tens_KE$.
\end{lemma}
\begin{demo}
  Since $K$ is complete and non-archimedean, $L\simeq K[T]/(f)$ for
  some Eisenstein polynomial $f\in\O_K[T]$, and this $f$ is still
  Eisenstein in $E[T]$ since $E/K$ is unramified. Therefore
  $L\tens_KE\simeq E[T]/(f)$ is a field, and $\O_E[T]/(f)$ is a DVR
  with fraction field $L\tens_KE$, whose maximal ideal is the image of
  $(T)$. This proves the first statement, and the second is clear
  since $L\tens_KE$ is a finite extension of $E$.
\end{demo}

\begin{lemma}\label{lemma:sigma-exists}
  There is a unique automorphism $\sigma$ of $K^\nr$ which induces the
  $q$th power Frobenius on $\bar k$ and whose fixed field is $K$.
\end{lemma}
\begin{demo}
  Suppose first that $K$ has characteristic $0$. The maximal
  absolutely unramified subextension of $K$ is isomorphic to
  $K_0=W(k)[p^{-1}]$, which we identify with a subfield of $K$. If we
  set $K_\infty=W(\bar k)[p^{-1}]$ there is a unique homomorphism
  $K_\infty\to K^\nr$ inducing the identity on residue fields. Lemma
  \ref{lemma:totally-ramified-extensions} shows that
  $K\tens_{K_0}K_\infty$ is a complete discretely valued field with
  residue field $\bar k$, so the canonical maps $K\to K^\nr$,
  $K_\infty\to K^\nr$ induce an isomorphism
  \begin{displaymath}
    K\tens_{K_0}K_\infty\isom K^\nr.
  \end{displaymath}
  We can then take $\sigma=1\tens F^f$, where $F$ is induced by the
  Frobenius of $W(\bar k)$ and $f$ satisfies $q=p^f$. This proves
  existence, and uniqueness follows from the isomorphism. For the last
  statement we let $e=[K:K_0]$ and pick a uniformizer $\pi$ of $K$, so
  that $1,\pi,\ldots,\pi^{e-1}$ is a basis of $K^\nr$ as $K_0$-vector
  space; since $\pi\in K$ is fixed by $\sigma$, the fixed field of
  $\sigma$ is exactly $K$.

  Suppose now $K$ has characteristic $p$ and again pick a uniformizer
  $\pi$ of $\O_K$. This time $k$ is a subfield of $\O_K$ in a
  canonical way, since $k$ is perfect, and $\pi$ is transcendental
  over $k$. Then $\O_K\simeq k[[\pi]]$ and from this we get
  $\O_{K^\nr}\simeq\bar k\ctens_k\O_K\simeq\bar k\ctens_kk[[\pi]]$
  where the topology on $k[[\pi]]$ is the $\pi$-adic one. The proof of
  existence and uniqueness is as before, and since
  $K^\nr\simeq\bar k((\pi))$ the last assertion is clear.
\end{demo}

From now on $\sigma$ will denote the automorphism of the lemma
\ref{lemma:sigma-exists} (it will also be used for several other
closely related maps). Suppose $L/K$ is a finite extension, and set
\begin{equation}
  \label{eq:L_K}
  L_K=K^\nr\tens_KL.
\end{equation}
We regard this as an $(K^\nr,L)$-bimodule via the natural maps
$K^\nr\to L_K$, $L\to L_K$.  The automorphism $\sigma$ of $K^\nr$
induces by functoriality an $L$-linear automorphism of $L_K$ which we
will also (as warned) denote by $\sigma$. The ring $L_K$ is a finite
direct sum of copies of the completion of a maximal unramified
extension of $L$; it will be worth the trouble to express this in an
intrinsic way. From now on $f$ denotes the residual degree of $L/K$.
Let $L_0\sset L$ be the maximal unramified extension of $K$ in $L$,
and denote by $\Sigma$ the set of $K$-embeddings $L_0\to K^\nr$; then
$[L_0:K]=|\Sigma|=f$. For $\iota\in\Sigma$ we denote by $K^\nr_\iota$
the field $K^\nr$ viewed as a $(K^\nr,L_0)$-bimodule where the
$L_0$-module structure is induced by $\iota$. We then set
\begin{equation}\label{eq:Lnr-from-Knr}
  L^\nr_\iota=K^\nr_\iota\tens_{L_0}L=K^\nr\tens_{\iota,L_0}L.
\end{equation}
Since $L/L_0$ is totally ramified, lemma
\ref{lemma:totally-ramified-extensions} shows that for all
$\iota\in\Sigma$, $L^\nr_\iota$ is a completion of a maximal unramfied
extension of $L$ (in particular the $L^\nr_\iota$ are isomorphic
extensions of $L$). It is also a $(K^\nr,L)$-bimodule compatible with
the $(K^\nr,L_0)$-bimodule structure of $K^\nr_\iota$. Now
\begin{equation}
  \label{eq:Knr-splitting}
  L_K\simeq (K^\nr\tens_KL_0)\tens_{L_0}L
\end{equation}
and 
\begin{align*}
  K^\nr\tens_KL_0&\isom\bigoplus_{\iota\in\Sigma}K^\nr_\iota\\
  (x,y)&\mapsto (x\iota(y))_{\iota\in\Sigma}
\end{align*}
so \ref{eq:Lnr-from-Knr} and \ref{eq:Knr-splitting} yield an
isomorphism
\begin{equation}
  \begin{split}
    \label{eq:L_K-splitting}
    L_K&\isom\bigoplus_{\iota\in\Sigma}L^\nr_\iota\\
    x\tens_Ky\tens_{L_0} b&\mapsto (x\iota(y)\tens_{L_0}b)_{\iota\in\Sigma}
  \end{split}
\end{equation}
where $x\in K^\nr$, $y\in L_0$ and $b\in L$ in the second line. This
splitting is compatible with the $(K^\nr,L)$-bimodule structures on
both sides.

For $\iota\in\Sigma$ the commutative diagram
\begin{displaymath}
  \xymatrix@R-0.5cm{
    K^\nr\ar[dd]_\sigma\\
    &L_0\ar[ul]_\iota\ar[dl]^{\sigma\iota}\ar[r]&L\\
    K^\nr
  }
\end{displaymath}
shows that that there is a homomorphism
\begin{displaymath}
  \sigma:L^\nr_\iota\to L^\nr_{\sigma\iota}
\end{displaymath}
(another use of $\sigma$!) making
\begin{displaymath}
  \xymatrix{
    K^\nr\ar[r]\ar[d]_\sigma&L^\nr_\iota\ar[d]^\sigma\\
    K^\nr\ar[r]&L^\nr_{\sigma\iota}
  }
\end{displaymath}
commutative. In terms of \ref{eq:Lnr-from-Knr} it is
\begin{displaymath}
  \sigma(x\tens_{\iota,L_0}b)=\sigma(x)\tens_{\sigma\iota,L_0}b.
\end{displaymath}
We can now express $\sigma:L_K\to L_K$ in terms
of the decomposition \ref{eq:L_K-splitting}, whose second line shows
that
\begin{equation}
  \label{eq:sigma-on-L_K-split}
  \sigma((x_\iota)_{\iota\in\Sigma})
  =(\sigma(x_{\sigma^{-1}\iota}))_{\iota\in\Sigma}.
\end{equation}
Recall that $f=[L_0:K]$, so the automorphism
\begin{equation}
  \label{eq:sigma-L}
  \sigma_L:=\sigma^f:L^\nr_\iota\to L^\nr_\iota
\end{equation}
is the identity of $L^\nr_\iota$ and induces the $q^f$-power
Frobenius on the residue field of $L^\nr_\iota$; i.e. it is the map
$\sigma$ of lemma \ref{lemma:sigma-exists} for $L^\nr_\iota/L$.

The field norm $N_{L/K}:L^\times\to K^\times$ extends to a homomorphism
\begin{equation}
  \label{eq:norm-LKK}
  N_{L/K}:L_K^\times\to (K^\nr)^\times
\end{equation}
which is the norm for the finite free $K^\nr$-algebra $K^\nr\to
L_K$. If $x\in L_K^\times$ corresponds $(x_\iota)$ under the decomposition
\ref{eq:L_K-splitting}, $N_{L/K}(x)$ is given by
\begin{equation}
  \label{eq:crys-norm}
  N_{L/K}(x)=\prod_{\iota\in\Sigma}N_{L^\nr_\iota/K^\nr}(x_\iota).
\end{equation}
We denote by $w_{L/K}:L_K^\times\to\bQ$ the homomorphism
\begin{equation}
  \label{eq:valuation-of-L_K}
  w_{L/K}(x)=[L:K]^{-1}v(N_{L/K}(x))
\end{equation}
where the numerical factor guarantees that $w_{L/K}$ extends the
valuation of $K^\times\subset L_K^\times$.  If $e$ is the ramification
index of $L/K$ then $[L^\nr_\iota:K^\nr]=e$ and $[L:K]=ef$, so
\ref{eq:valuation-of-L_K} can be written
\begin{equation}
  \label{eq:w_L/K}
  w_{L/K}(x)=\sum_{\iota\in\Sigma}f^{-1}v(x_\iota).
\end{equation}

\subsection{Lubin-Tate $F$-isocrystals.}
\label{sec:lubin-tate-fisoc}

From now on $L/K$ will be a finite \emph{Galois} extension with degree
$d$ and group $G$ (for which we will write $G_{L/K}$ if other Galois
groups are in the picture). We fix a uniformizer $\pi$ of $\O_K$; it
is fixed by $\sigma$, so we can now apply
\S\ref{sec:F-isocrystals}. It will be convenient to make a slight
change of notation: we will write $V_K(\lambda)$ and $D_K(\lambda)$
for what were previously denoted $V_{K^\nr}(\lambda)$ and
$D_{K^\nr}(\lambda)$. We will sometimes drop the subscript $K$ if the
meaning is clear.

By definition the fundamental class
$u_{L/K}\in H^2(G,L^\times)$ is an element of invariant $1/d$; in
other words it is the class of a central simple $K$-algebra of
invariant $1/d$. By theorem \ref{thm:end-ring-of-simple-Fisocrystal}
we can take this $K$-algebra to be the endomorphism algebra of any
$F$-isocrystal $V$ on $K^\nr$ of rank $d$ and invariant $-1/d$. Since
we want to construct a $G$-valued 2-cocycle for $u_{L/K}$ we will need
an explicit splitting of $\End(V)$ over $L$, which amounts to giving
and embedding $L\to\End(V)$ over $K$.

We first note that the action of $G$ on $L$ induces an action on
$L_K=K^\nr\tens_KL$, and this action commutes with $\sigma$. If
$\alpha\in L_K^\times$,
\begin{equation}
  \label{eq:our-basic-Frobenius}
  F_\alpha:L_K\to L_K\qquad
  F_{\alpha}(x)=\alpha\cdot\up\sigma x 
\end{equation}
defines a Frobenius structure on $L_K$, and we denote by $L_K(\alpha)$
the $F$-isocrystal $(L_K,F_\alpha)$. We denote by $D_K(\alpha)$ the
endomorphism ring of $L_K(\alpha)$. Since $\sigma\tens1$ is
$L$-linear, right multiplication defines an embedding
$L\to D_K(\alpha)$.

\begin{prop}\label{prop:slope-of-M(alpha)}
  For any finite Galois extension $L/K$ and $\alpha\in L_K^\times$,
  the $F$-isocrystal $L_K(\alpha)$ is isopentic of slope
  $w_{L/K}(\alpha)$.
\end{prop}
\begin{demo}
  Write $\alpha=(a_\iota)$; if $x=(x_\iota)$ then
  \begin{displaymath}
    F_\alpha(x)
    =\alpha\sigma(x)=(a_\iota\sigma(x_{\sigma^{-1}\iota})
  \end{displaymath}
  and a simple induction yields
  \begin{displaymath}
    F_\alpha^m(x)=(a_\iota\sigma(a_{\sigma^{-1}\iota})\cdots
    \sigma^{m-1}(a_{\sigma^{-(m-1)}\iota})\sigma^m(x_{\sigma^{-m}\iota})).
  \end{displaymath}
  In particular since $\sigma^f\iota=\iota$ and $\sigma^f=\sigma_L$,
  \begin{displaymath}
    F_\alpha^f=(b_\iota\sigma_L(x_\iota))
  \end{displaymath}
  where
  \begin{displaymath}
    b_\iota=a_\iota\sigma(a_{\sigma^{-1}\iota})\cdots
    \sigma^{f-1}(a_{\sigma^{-(f-1)}\iota}).
  \end{displaymath}
  Since $\Sigma=\{\sigma^{-i}\iota\,|\, 0\le i<f\}$ and
  $\sigma:L^\nr_\iota\to L^\nr_{\sigma\iota}$ preserves valuations,
  \ref{eq:w_L/K} shows that
  \begin{displaymath}
    v_L(b_\iota)=\sum_{0\le i<f}v(a_{\sigma^{-i}\iota})
    =fw_{L/K}(\alpha).
  \end{displaymath}
  Let $e$ be the ramification index of $L/K$; then $d=ef$,
  \begin{displaymath}
    F^d(x)=F^e(b_\iota\cdot\up{\sigma_L}x_\iota)
    =(\up{1+\sigma_L+\cdots+\sigma_L^{e-1}}b_\iota\up{\sigma_L^e}x_\iota)
  \end{displaymath}
  and
  \begin{displaymath}
    v(\up{1+\sigma_L+\cdots+\sigma_L^{e-1}}b_\iota)
    =ev(b_\iota)=efw_{L/K}(\alpha)
    =dw_{L/K}(\alpha).
  \end{displaymath}

  Finally let $M$ be the lattice
  \begin{displaymath}
    M=\bigoplus_\iota\O_{L^\nr_\iota}\subset\bigoplus_\iota L^\nr_\iota=L_K.
  \end{displaymath}
  Since $dw_{L/K}(\alpha)\in\bZ$, the last few calculations show that
  \begin{displaymath}
    F^d(M)=\pi^{dw_{L/K}(\alpha)}M
  \end{displaymath}
  for any uniformizer $\pi$ of $K$. Then
  $F^{nd}(M)=\pi^{ndw_{L/K}(\alpha)}M$, and proposition
  \ref{prop:bounds-on-Frobenius} shows that $(L_K,F_\alpha^d)$ has
  slope $dw_{L/K}(\alpha)$. By the remark at the end of section
  \ref{sec:F-isocrystals}, $(L_K,F_\alpha)$ has slope
  $w_{L/K}(\alpha)$, as asserted.
\end{demo}

\subsubsection{The unramified case.}
\label{sec:LK(lambda)-unramified}

Let's consider the case when $L/K$ unramified of degree $d=f$. Pick a
$r\in\bZ$ that is relatively prime to $d$ and set
\begin{equation}
  \label{eq:unramified-alpha}
  \alpha=(\pi^r,1,\ldots,1)
\end{equation}
where $\pi\in K$ is a uniformizer and we use the identification
\ref{eq:L_K-splitting}. By the proposition $L_K(\alpha)$ has slope
$r/d$, and as it has rank $d$ it must be isomorphic to $V_K(r/d)$. In
fact this is easy to write down an explicit isomorphism from the
explicit formula for $F_\alpha$ in the proof of proposition
\ref{prop:slope-of-M(alpha)}.

\subsection{The fundamental extension.}
\label{sec:extension-class}

We can now identify $D_K(\alpha)$ with a crossed product algebra
associated to a certain an $L^\times$-valued $2$-cocycle for $G$. In
the case of the fundamental class we then identify this $2$-cocycle as
the class of a certain Yoneda $2$-extension.

\begin{lemma}\label{lemma:solving-1-F}
  The sequence
  \begin{displaymath}
    1\to L^\times\to L^\times_K\Xto{\sigma-1}L^\times_K
    \Xto{w_{L/K}}d^{-1}\bZ\to0
  \end{displaymath}
  is exact.
\end{lemma}
\begin{demo}
  The only non-obvious point is that the image of $\sigma-1$ is the
  kernel of $w_{L/K}$, and it is clearly contained in the kernel. We
  first observe that the subgroup of $\bZ^\Sigma$ consisting of
  $(n_\iota)$ such that $\sum_\iota n_\iota=0$ is spanned by elements
  with a $1$ in some position $\iota$, $-1$ in position
  $\sigma^{-1}\iota$ and zeros elsewhere. From this it follows that
  any $(x_\iota)\in L_K^\times$ in the kernel of $w_{L/K}$ is
  congruent modulo the image of $\sigma-1$ to a $(x_\iota)$ such that
  $v_L(x_\iota)=0$ for all $\iota$. Suppose $x=(x_\iota)$ is one such
  element of $L_K$; since $\sigma_L-1$ is a retraction of
  $(L^\nr_\iota)^\times$ onto its subgroup of elements of valuation
  zero, for each $\iota$ there is a $c_\iota$ such that
  $x_\iota=\sigma_L(c_\iota)/c_\iota$. Then $x=\sigma_L(c)/c$ in
  $L_K^\times$ where $c=(c_\iota)$, and since
  \begin{displaymath}
    \sigma_L-1=\sigma^f-1=(\sigma-1)(\sigma^{f-1}+\cdots+\sigma+1)
  \end{displaymath}
  we conclude that $x$ is in the image of $\sigma-1$.
\end{demo}

Recall that $L/K$ is Galois with group $G$; since
$w_{L/K}(\up{s-1}\alpha)=0$, the lemma shows that for any $s\in G$,
there is a $\beta_s\in L_K^\times$ such that
\begin{equation}
  \label{eq:beta-to-alpha}
  \up{\sigma-1}\beta_s=\up{s-1}\alpha
\end{equation}
well defined up to a factor in $L^\times$. From \ref{eq:beta-to-alpha}
we see that
\begin{equation}
  \label{eq:u-s}
  u_s:L_K\to L_K\qquad
  u_s(x)=\beta_s\cdot\up sx
\end{equation}
commutes with $F_\alpha$:
\begin{align*}
  F_\alpha(u_s(x))&=\alpha\cdot\up\sigma(\beta_s\cdot\up sx)\\
  &=\alpha\cdot\up\sigma\beta_s\cdot\up{\sigma s}x
  =\up s\alpha\cdot\beta_s\cdot\up{s\sigma}x\\
  &=\beta_s\cdot\up s(\alpha\cdot\up\sigma x)=u_s(F_\alpha(x)).
\end{align*}
and is thus an automorphism of $L_K(\alpha)$.

A quick calculation shows that the $u_s$ satisfy the relations
\begin{equation}
  \label{eq:crossed-product}
  u_su_t=a_{s,t}u_{st},\quad
  u_s\ell=\up s\ell u_s\quad
\end{equation}
for all $\ell\in L$ and $s$, $t\in G$, where 
\begin{equation}
  \label{eq:2-cocycle}
  a_{s,t}=\up s\beta_t\beta_{st}^{-1}\beta_s.
\end{equation}
Since $\sigma$ commutes with the action of $G$, equation
\ref{eq:beta-to-alpha} shows that $\up{\sigma-1}a_{s,t}=1$ and thus
$a_{s,t}\in L^\times$ by lemma \ref{lemma:solving-1-F}. Thus
$(a_{s,t})$ is an $L^\times$-valued 2-cocyle for $G$; we take can
$\beta_1=1$, making $(a_{s,t})$ a normalized cocycle, which we do from
now on.

\begin{thm}\label{thm:the-fundamental-class}
  The endomorphism algebra $D_{L/K}(\alpha)$ is isomorphic to the
  crossed product algebra associated to the 2-cocycle
  \ref{eq:2-cocycle} associated to $\alpha$. The extension
  corresponding to $D_{L/K}(\alpha)$ is isomorphic to
  \begin{displaymath}
    1\to L^\times\Xto{i}W(\alpha)\Xto{p}G\to1
  \end{displaymath}
  where $W(\alpha)$ is the set of pairs
  $(\beta,s)\in L_K^\times\times G$ satisfying
  \begin{equation}
    \label{eq:elements-of-W}
    \up{\sigma-1}\beta=\up{s-1}\alpha
  \end{equation} with the
  composition law
  \begin{equation}
    \label{eq:group-law-of-W}
    (\beta,s)(\gamma,t)=(\beta\cdot\up s\gamma,st)
  \end{equation}
  and $i$, $p$ are the evident inclusion and projection.
\end{thm}
\begin{demo}
  The \ref{eq:crossed-product} are the defining relations for the
  crossed product algebra $A$ defined by the cocycle
  \ref{eq:2-cocycle}. The universal property of such algebras (e.g.\
  \cite[\S16 $n^o$ 9 Prop. 12]{bourbaki:csa}) yields a nonzero
  $K$-homorphism $A\to D_K(\alpha)$ (since nonzero on $K$) which is
  necessarily an isomorphism since both sides have dimension $d^2$.

  For the second part we need only remark that $i:K\to D_K(\alpha)$ is
  the canonical injection, the extension group corresponding to
  $D_K(\alpha)$ is the set of $u\in D_{L/K}(\alpha)$ satisfying
  $i\circ s=\ad(u)\circ i$ for some $s\in G$, and it is easily checked
  that these have the form $u=(\beta,s)$ with $\beta\in L_K^\times$
  satisfying \ref{eq:elements-of-W}.
\end{demo}

\begin{cor}\label{cor:the-fundamental-class}
  For any finite Galois extension $L/K$ of degree $d$ and any
  $\alpha\in L^\times_K$ such that $w_{L/K}(\alpha)=-1/d$, the
  2-cocycle \ref{eq:2-cocycle} represents the fundamental class of
  $L/K$.
\end{cor}

When $w_{L/K}(\alpha)=-1/d$, the extension $W(\alpha)$ representing the
fundamental class is by definition the \textit{Weil group} of $L/K$
and denoted by $W_{L/K}$. Up to isomorphism of extensions it is
independent of the choice of $\alpha$; in fact since
$H^1(G,L^\times)=0$ it is well-defined up to inner automorphisms by
elements of $L^\times$.


We note finally that just as the Galois group $\Gal(L/K)$ is (by
definition!) an automorphism group, the same is true for $W(\alpha)$:
by construction it consists of automorphisms $u$ of $L_K(\alpha)$
satisfying $i\circ s=\ad(u)\circ i$ for some $s\in G$. When $K=\bQ_p$,
this observation is due to Morava \cite{morava:1979}.

\subsubsection{The unramified case again.}
\label{sec:unramified-cocycle}

Suppose $L/K$ is unramified, so that $f=d$ and we identify
$K^\nr\simeq L^\nr$. The Galois group $G=\Gal(L/K)$ is generated by
the ``arithmetic Frobenius'' $\sigma_\arith=1\tens(\sigma|L)$ and we
identify $G\simeq\bZ/d\bZ$ by means of the generator $\sigma_\arith$;
in the formulas to follow we also identify elements of $\bZ/d\bZ$ with
integers in the range $[0,d)$. Since $L/K$ is unramified we can take
$\alpha=(a^{-1},1,\ldots,1)$ with $a\in K^\times$ (e.g. $a=\pi^r$ as
in equation \ref{eq:unramified-alpha}) Then
\ref{eq:sigma-on-L_K-split} yields
\begin{displaymath}
  \up{\sigma_\arith^i-1}\alpha=(a,1,\ldots,1,
  \underset{\underset{d-i}\uparrow}a^{-1},1,\ldots,1)
\end{displaymath}
and a $\beta_i\in L_K^\times$
satisfying $\up{\sigma-1}\beta_i=\up{\sigma_\arith^i-1}\alpha$ is
\begin{equation}
  \label{eq:unramified-beta}
  \beta_i=(a^{-1},a^{-1},\ldots,\!
  \underset{\underset{d-i-1}\uparrow}{a^{-1}}\!\!,
  1,\ldots,1).
\end{equation}
From this it follows easily that
\begin{displaymath}
  \beta_i(\sigma_\arith^i(\beta_j))\beta_{i+j}^{-1}=
  \begin{cases}
    a^{-1}&i+j<d\\ 1&i+j\ge d
  \end{cases}.
\end{displaymath}
To normalize this cocycle we note that the coboundary of the constant
1-cochain with value $a$ is the constant 2-cocycle with value
$a$. Adding this to the previous cocycle yields the well-known
2-cocycle
\begin{displaymath}
  a_{i,j}=
  \begin{cases}
    1&i+j<d\\ a&i+j\ge d
  \end{cases}
\end{displaymath}
representing the fundamental class in the unramified case.

A class in $\Br{L/K}$ gives rise, via the isomorphisms
\begin{equation}
  \label{eq:Br-as-H2}
  \Br{L/K}\simeq H^2(G,L^\times)\simeq\Ext^2_G(\bZ,L^\times)
\end{equation}
to a Yoneda 2-extension of $\bZ$ by $L^\times$. We can identify the
extension corresponding to the fundamental class:

\begin{thm}\label{thm:2-extension}
  Suppose $[L:K]=d$ and write $w=w_{L/K}$. The class of
  $u_{L/K}\in\Br{L/K}$ is the opposite of the class of the extension
  \begin{equation}
    \label{eq:2-extension}
    0\to L^\times\to
    L_K^\times\Xto{\sigma-1}L_K^\times
    \Xto{dw}\bZ\to 0.
  \end{equation}
  in $\Ext^2_G(\bZ,L^\times)$.
\end{thm}
\begin{demo}
  By lemma \ref{lemma:solving-1-F}, \ref{eq:2-extension} is exact and
  we just need to recall the usual recipe (e.g.\ \cite[\S7 $n^o$
  3]{bourbaki:hom}) for associating to it a class in
  $\Ext^2_G(\bZ,L^\times)\simeq H^2(G,L^\times)$. If $P_\cdot\to\bZ$
  is a resolution of $\bZ$ by projective $\bZ[G]$-modules, there is a
  morphism of complexes
  \begin{equation}
    \label{eq:calculation-of-extension-class}
    \xymatrix{
      \ar[r]
      &P_3\ar[r]^{d_3}\ar[d]^{f_3}
      &P_2\ar[r]^{d_2}\ar[d]^{f_2}&P_1\ar[r]^{d_1}\ar[d]^{f_1}
      &P_0\ar[r]^{d_0}\ar[d]^{f_0}&\bZ\ar[r]\ar@{=}[d]&0\\
      \ar[r]
      &1\ar[r]
      &L^\times\ar[r]&L_K^\times\ar[r]^{\sigma-1}
      &L_K^\times\ar[r]^{dw}
      &\bZ\ar[r]&0
    }
  \end{equation}
  unique up to homotopy. The map $f_2$ lies in the kernel of
  $\Hom(P_2,L^\times)\to\Hom(P_3,L^\times)$, and its image in
  $H^2(G,L^\times)$ is the class of \ref{eq:2-extension}. For
  $P_\cdot$ we take the standard bar complex:
  \begin{displaymath}
    d_0([\ ])=1,\qquad
    d_1([s])=s[\ ]-[\ ],\qquad
    d_2([s|t])=s[t]-[st]+[s].
  \end{displaymath}
  If we set
  \begin{displaymath}
    f_0([\ ])=\alpha,\qquad
    f_1([s])=\beta_s,\qquad
    f_2([s|t])=a_{s,t}
  \end{displaymath}
  the commutativity of \ref{eq:calculation-of-extension-class}
  shows that 
  \begin{displaymath}
    dw(\alpha)=1,\qquad
    \up{\sigma-1}\beta_s=\up{s-1}\alpha,\qquad
    a_{s,t}=\up s\beta_t\beta_{st}^{-1}\beta_s.
  \end{displaymath}
  The first equality shows that $D_{L/K}(\alpha)$ represents
  $-u_{L/K}$, and the others show that $(a_{s,t})$ represents
  $D_{L/K}(\alpha)$ in $H^2(G,L^\times)$. Since the class of
  \ref{eq:2-extension} is $(a_{s,t})$ we are done.
\end{demo}

\subsection{The norm residue symbol.}
\label{sec:norm-residue-symbol}

We denote by
\begin{equation}
  \label{eq:inverse-reciprocity}
  \eta_{L/K}:G^\ab\to K^\times/N_{L/K}(L^\times)
\end{equation}
the map induced by the cup product
\begin{displaymath}
  \cup u_{L/K}:\hat H^{-2}(G,\bZ)\to\hat H^0(G,L^\times)
\end{displaymath}
and the standard idenfication of these Tate cohomology groups with
those in \ref{eq:inverse-reciprocity}. We know that $\eta_{L/K}$ is an
isomorphism and we denote its inverse, the norm residue symbol by
\begin{equation}
  \label{eq:norm-residue-symbol}
  \theta_{L/K}=\eta_{L/K}^{-1}:K^\times/N_{L/K}(L^\times)\to G^\ab.
\end{equation}
We will also write $\theta_{L/K}$ for the composite
\begin{displaymath}
  K^\times\to K^\times/N_{L/K}L^\times\to G^\ab.
\end{displaymath}
and we will also use the traditional notation for
\begin{equation}
  \label{eq:norm-residue-symbol-traditional}
  (a,L/K)=\theta_{L/K}(a)
\end{equation}
where $a$ is an element of either $K^\times$ or
$K^\times/N_{L/K}L^\times$.

That $\eta_{L/K}$ is an isomorphism is of course a consequence of the
Tate-Nakayama theorem, but in fact this follows from the exact
sequence \ref{eq:2-extension} and standard arguments (as in
\cite[Ch. XIII \S5 Ex. 2]{serre:1968}) once it is shown that
$L_K^\times$ is a cohomologically trivial $G$-module. In fact if
$I\sset G$ is the inertia subgroup, $L_K^\times$ is the $G$-module
induced from the $I$-module $(L^\nr)^\times$, and the latter is
cohomologically trivial since the norm
$(L^\nr)^\times\to(K^\nr)^\times$ is surjective. We should point out
that the role of the extension \ref{eq:2-extension} is obscured in
\cite{serre:1968} since in the case that Serre is considering
(quasi-finite residue fields) the sequence
$L_K^\times\Xto{1-F}L_K^\times\to\bZ$ is \textit{not} always
exact. However the argument works since the cohomology group here is
cohomologically trivial; c.f.\ \cite[\textit{loc.\ cit.}]{serre:1968}.

\begin{thm}\label{thm:dwork-serre-formula}
  Suppose $L/K$ is a Galois extension of local fields with group $G$
  and set $d=[L:K]$. Let $\alpha$ be an element of $L_K^\times$ such
  that $w_{L/K}(\alpha)=-1/d$, and let $s\in G$. If
  $\beta\in L_K^\times$ satisfies
  \begin{displaymath}
    \up{\sigma-1}\beta=\up{s-1}\alpha    
  \end{displaymath}
  then
  \begin{equation}
    \label{eq:dwork-serre-formula}
    \eta_{L/K}(s)=N_{L/K}(\beta)\ \mod N_{L/K}L^\times.
  \end{equation}
\end{thm}
\begin{demo}
  This follows from the formula \cite{serre:1968} for the cup product
  and the definition \ref{eq:2-cocycle}:
  \begin{displaymath}
    \bar s\cup u_{L/K}=\prod_{t\in G}a_{t,s}
    =\prod_{t\in G}\up t\beta_s\beta_{ts}^{-1}\beta_t=N_{L/K}(\beta_s)
  \end{displaymath}
  since we can take $\beta=\beta_s$.
\end{demo}

When $L/K$ is totally ramified we may identify $L_K=L^\nr$ and take
$\alpha=\pi^{-1}$ for any uniformizer $\pi$ of $L^\nr$; in this case
\ref{eq:dwork-serre-formula} is equivalent to the formula proven by
Dwork \cite{dwork:1958}.

\subsection{Successive extensions.}
\label{sec:reciprocity-successive-extensions}

The formal properties of the norm residue symbol relative to
successive extensions follow from theorem
\ref{thm:dwork-serre-formula} by direct computations. This could have
been left as an exercise, but the constructions we need for this will
also be used in section \ref{sec:shafarevich-weil}.

We first extend the constructions of section \ref{sec:L_K}.  Suppose
$L/K$ and $E/L$ are finite extensions of local fields. The
transitivity isomorphism $K^\nr\tens_KE\simeq(K^\nr\tens_KL)\tens_LE$
may be written
\begin{equation}
  \label{eq:transitivity}
  E_K\simeq L_K\tens_LE.
\end{equation}
This shows that $E_K$ has a canonical structure of a free
$L_K$-module, and we denote by
\begin{equation}
  \label{eq:crys-norm3}
  N_{E/L}^K:E_K^\times\to L_K^\times  
\end{equation}
the norm map. When $E/L/K$ is $L/K/K$ this is the norm
map \ref{eq:norm-LKK} introduced earlier. 

If $F/E$ is another finite extension, the norm is transitive:
\begin{equation}
  \label{eq:norm-transitive}
  N^K_{E/L}\circ N^K_{F/E}=N^K_{F/L}.
\end{equation}
Applying this in the case $E/L/K/K$ and invoking
\ref{eq:valuation-of-L_K}, we find that
\begin{equation}
  \label{eq:norm-valuation-compatibility}
  w_{L/K}(N_{E/L}^K(\beta))=[E:L]w_{E/K}(\beta)
\end{equation}
for $\beta\in E_K$.

The norm is equivariant for isomorphisms in the sense that
if $f:E\to E'$ is an isomorphism of local fields and $L'=f(L)$,
$K'=f(K)$ then
\begin{equation}
  \label{eq:norm-and-galois}
  \xymatrix@C+1em{
    E_K^\times\ar[r]^{N^K_{E/L}}\ar[d]_f&L_K^\times\ar[d]^f\\
    (E')_{K'}^\times\ar[r]_{N^{K'}_{E'/L'}}&(L')_{K'}^\times
  }
\end{equation}
is commutative.

For the rest of this section $E/K$ is a finite Galois extension, and
$K\sset L\sset E$. Our notation for Galois groups will now indicate
the extension, so for example the Galois group of $E/K$ is $G_{E/K}$.

We can now deduce two of the four basic properties of the norm residue
symbol relative to this situation.  Suppose that $L/K$ is Galois and
denote by
\begin{displaymath}
  \pi^K_{E/L}:G_{E/K}\to G_{L/K}
\end{displaymath}
the canonical homomorphism. To calculate $\eta_{E/K}$ we choose
$\alpha\in E_K^\times$ such that $w_{E/K}(\alpha)=1/[E:K]$: if $s\in
G_{E/K}$ we pick 
$\beta\in E_K^\times$ such that
$\up{\sigma-1}\beta_s=\up{s-1}\alpha$; then
$\eta_{E/K}(s)=N_{L/K}(\beta_s)$ (note that $N_{L/K}=N^K_{L/K}$
here). On
the other hand \ref{eq:norm-valuation-compatibility} says that
\begin{displaymath}
   w_{L/K}(N^K_{E/L}(\alpha))=\frac{[E:L]}{[E:K]}=\frac{1}{[L:K]}
\end{displaymath}
and we may use $N_{E/L}(\alpha)$ to calculate $\eta_{L/K}$. 
By equivariance
\begin{equation}
  \label{eq:eta-quotient}
  \up{\sigma-1}N_{E/L}(\beta)=\up{s-1}N_{E/L}(\alpha)
\end{equation}
so from \ref{eq:norm-transitive} and \ref{eq:dwork-serre-formula} we
get
\begin{equation}
  \label{eq:eta-quotient1}
  \eta_{E/K}(s)=\eta_{L/K}(\pi^K_{E/L}(s))\mod N_{L/K}L^\times
\end{equation}
which is equivalent to the formula
\begin{equation}
  \label{eq:norm-residue-quotients}
  \pi^K_{E/L}(a,E/K)=(a,L/K).
\end{equation}

The functoriality of the Artin symbol follows from our constructions
and \ref{eq:norm-and-galois} applied with $E'=E$ and $f=t\in
G_{E/K}$. Given $s$, $t\in G_{E/K}$ and
$\up{\sigma-1}\beta_s=\up{s-1}\alpha$, the equality
\begin{displaymath}
  \up{\sigma-1}(\up t\beta_s)=\up t(\up{\sigma-1}\beta_s)
  =\up{ts-t}\alpha=\up{tst^{-1}-1}(\up t\alpha)
\end{displaymath}
implies that
\begin{equation}
  \label{eq:eta-functorial1}
  \eta_{E/t(L)}(tst^{-1}))=\up t\eta_{E/L}(s)
\end{equation}
or
\begin{equation}
  \label{eq:eta-functorial2}
  t(a,E/t(L))t^{-1}=(\up ta,E/L)
\end{equation}
modulo appropriate subgroups.

To figure out the behavior of the Artin symbol with respect to the
inclusion $G_{E/L}\to G_{E/K}$ we need to understand how $E_K$ and
$E_L$ are related. In fact the same construction used to decompose
$L_K$ into a sum of isomorphic copies of $L^\nr$ also shows that $E_K$
is a direct sum of copies of $E_L$, and in fact reduces to it when
$E=L$.

Recall that for any $K$-embedding $\iota:L_0\to K^\nr$ we defined
$L_\iota^\nr=K^\nr_\iota\tens_{L_0}L$, and we now set
\begin{displaymath}
  E_{\iota,L}=L^\nr_\iota\tens_LE\simeq(K^\nr_\iota\tens_{L_0}L)\tens_LE
  \simeq K^\nr_\iota\tens_{L_0}E.
\end{displaymath}
The $E_{\iota,L}$ are isomorphic as $E$-algebras. On the other hand
\begin{equation}\label{eq:EK-and-EL}
  \begin{split}
    E_K&=K^\nr\tens_KE\simeq(K^\nr\tens_KL_0)\tens_{L_0}E\\
       &\simeq\bigoplus_{\iota\in\Sigma}K^\nr_\iota\tens_{L_0}E\\
       &\simeq\bigoplus_{\iota\in\Sigma}E_{\iota,L}.
  \end{split}
\end{equation}
and this isomorphism is equivariant for the action of $G_{E/L}$ on
both sides. The $L$-linear isomorphisms
$\sigma:L^\nr_\iota\isom L^\nr_{\sigma\iota}$ induce $E$-linear isomorphisms
\begin{displaymath}
  \sigma:E_{\iota,L}\isom E_{\sigma\iota,L}
\end{displaymath}
and the formula for $\sigma:E_K\to E_K$ is just
\ref{eq:sigma-on-L_K-split} but where now $x_\iota\in
E_{\iota,L}$. Finally if we are given an extension $F/E$ the
preceeding applies to the $F/E/L$, so that
$F_K=K^\nr\tens_KF$ has a decomposition
\begin{equation}
  \label{eq:FK-decomp}
  F_K\simeq\bigoplus_{\iota\in\Sigma}F_{\iota,L}
\end{equation}
with $F_{\iota,L}$ is defined like $E_{\iota,L}$, and
$F_{\iota,L}\simeq E_{\iota,L}\tens_LF$. This leads to norm maps
\begin{displaymath}
  N^{L\,\iota}_{F/E}:F^\times_{\iota,L}\to E^\times_{\iota,L}
\end{displaymath}
and if $x\in F_K$ is identified with $(x_\iota)$ via
\ref{eq:FK-decomp}, 
\begin{equation}
  \label{eq:norm-in-3-fold-ext}
  N^L_{F/E}(x)=(N^{L,\iota}_{F/E}(x_\iota))_{\iota\in\Sigma}.
\end{equation}

For the next constructions we choose a ``base point''
$\iota_0\in\Sigma$ and identify $E_L\simeq E_{\iota_0,L}$. Then 
\begin{equation}
  \label{eq:gamma-ELK}
  \gamma^E_{L/K}:E_L^\times\to E_K^\times
\end{equation}
denotes the inclusion of the direct summand $E^\times_{\iota_0,L}$ in
the multiplicative version
\begin{displaymath}
  E_K^\times\simeq\bigoplus_{\iota\in\Sigma}E^\times_{\iota,L}.
\end{displaymath}
of \ref{eq:EK-and-EL}. Since $\sigma^f|L_0=1$, the distinct
$K$-embeddings $L_0\to K^\nr$ are the $\sigma^i\iota_0$ with
$i\in\bZ/f\bZ$. The ``twisted diagonal''
\begin{equation}
  \label{eq:delta-ELK}
  \delta^E_{L/K}:E_L^\times\to E_K^\times
\end{equation}
is defined as follows: for $i\in\bZ/f\bZ$ and $x\in E_L$, the
$\sigma^i\iota_0$-component of $\delta^E_{L/K}(x)$ is
$\sigma^i(x)$. Since we can identify
$\Sigma=\{\sigma^i\iota_0\,|\,i\in\bZ/f\bZ\}$ we can write the
decomposition \ref{eq:EK-and-EL} as indexed over $i\in\bZ/f\bZ$; in
this notation
\begin{equation}
  \label{eq:gamma-delta-maps}
  \begin{split}
    \gamma^E_{L/K}(x)&=(x,1,\ldots,1)\\
    \delta^E_{L/K}(x)&=(x,\sigma(x),\sigma^2(x),\ldots,\sigma^{f-1}(x))
  \end{split}
\end{equation}
and $\sigma:E_K\to E_K$ is
\begin{equation}
  \label{eq:sigma-in-Sigma}
  \sigma(x_0,\ldots,x_{f-1})=(\sigma(x_{f-1}),\sigma(x_0),\sigma(x_1),\ldots).
\end{equation}

Since $G_{E/L}$ fixes $L_0\sset L$ the decomposition
\ref{eq:EK-and-EL} is equivariant for the action $G_{E/L}$, and it
follows that $\gamma^E_{L/K}$ and $\delta^E_{L/K}$ are equivarient as
well. Corresponding to the formula
\ref{eq:norm-valuation-compatibility} we have, for
$\alpha\in E_L^\times$,
\begin{equation}
  \label{eq:iota-and-w}
  w_{E/K}(\gamma^E_{L/K}(\alpha))=\frac{1}{[L:K]}w_{E/L}(\alpha).
\end{equation}

For any 3-fold extension $F/E/L/K$ the maps $\gamma$ and $\delta$ are
compatible with the norms. In fact we see from
\ref{eq:norm-in-3-fold-ext} and \ref{eq:gamma-delta-maps} that
\begin{equation}
  \label{eq:iota-delta-norms}
  \xymatrix{
    F^\times_L\ar[r]^{\delta^F_{L/K}}\ar[d]_{N^L_{F/E}}&F^\times_K\ar[d]^{N^K_{F/E}}\\
    E^\times_L\ar[r]_{\delta^E_{L/K}}&E_K^\times
  }
  \qquad
  \xymatrix{
    F^\times_L\ar[r]^{\gamma^F_{L/K}}\ar[d]_{N^L_{F/E}}&F^\times_K\ar[d]^{N^K_{F/E}}\\
    E^\times_L\ar[r]_{\gamma^E_{L/K}}&E^\times_K
  }
\end{equation}
are commutative.

We need one more fact about these maps. The composite
\begin{displaymath}
  L_L^\times\Xto{\delta^L_{L/K}}L_K^\times\Xto{N^K_{L/K}}K_K^\times
\end{displaymath}
is a homomorphism $(L^\nr)^\times\to(K^\nr)^\times$. It is not the
usual norm map for $L^\nr/K^\nr$, however:

\begin{lemma}\label{lemma:random-norm-factoid}
  For $\ell\in L^\times\subset(L^\nr)^\times$,
  \begin{equation}
    \label{eq:random-norm-factoid}
    N^K_{L/K}(\delta^L_{L/K}(\ell))=N_{L/K}(\ell)
  \end{equation}
  where the norm on the right is the norm for $L/K$.
\end{lemma}
\begin{demo}
  From \ref{eq:gamma-delta-maps} and \ref{eq:crys-norm} we get
  \begin{displaymath}
    N^K_{L/K}(\delta^L_{L/K}(\ell))=\prod_{0\le i<f}\sigma^i(N_{L^\nr/K^\nr}(\ell)).
  \end{displaymath}
  Then $L^\nr\simeq L\tens_{L_0}K^\nr$ shows that
  $N_{L^\nr/K^\nr}(\ell)=N_{L/L_0}(\ell)$, so the product on the right
  hand side is
  \begin{displaymath}
    \prod_{0\le
      i<f}\sigma^i(N_{L/L_0}(\ell))=N_{L_0/K}(N_{L/L_0}(\ell))
    =N_{L/K}(\ell)
  \end{displaymath}
  since $\sigma|L_0$ generates the Galois group of $L_0/K$.
\end{demo}

We can now show that the inclusion $G_{E/L}\to G_{E/K}$ corresponds
via the Artin symbol to the map
$L^\times/N_{E/L}E^\times\to K^\times/N_{E/K}E^\times$ induced by the
norm. Take $s\in G_{E/L}$; to compute $\eta_{E/L}(s)$ we choose
$\alpha\in E_L^\times$ such that $w_{E/L}(\alpha)=[E:L]^{-1}$ and
$\beta\in E_L^\times$ such that $\up{\sigma_L-1}\beta=\up{s-1}\alpha$;
then $\eta_{E/L}(s)$ is the class of $N_{E/L}(\beta)$ in
$L^\times/N_{E/L}E^\times$ (recall here that $N_{E/L}$ here is
$N^L_{E/L}$). On the other hand by \ref{eq:iota-and-w} we have
$w_{E/K}(\gamma^E_{L/K}(\alpha))=1/[E:K]$, so we can use
$\gamma^E_{L/K}(\alpha)$ to compute $\eta_{E/K}(s)$. From
\ref{eq:gamma-delta-maps} and \ref{eq:sigma-in-Sigma} we get
\begin{equation}
  \label{eq:artin-symbol-subgroup}
  \up{\sigma-1}\delta^E_{L/K}(\beta)=(\up{\sigma_L-1}\beta,1,\ldots,1)
  =(\up{s-1}\alpha,1,\ldots,1)=\up{s-1}\gamma^E_{L/K}(\alpha)
\end{equation}
and thus
\begin{displaymath}
  \eta_{E/K}(s)=N^K_{E/K}(\delta^E_{L/K}(\beta)) \mod N_{E/K}E^\times.
\end{displaymath}
To evaluate the right hand side we consider the diagram
\begin{displaymath}
  \xymatrix{
    E_L^\times\ar[r]^{\delta^E_{L/K}}\ar[d]_{N^L_{E/L}}
    &E_K^\times\ar[d]_{N^K_{E/L}}\ar@/^0.5cm/[rd]^{N^K_{E/K}}\\
    L_L^\times\ar[r]_{\delta^L_{L/K}}&L_K^\times\ar[r]_{N^K_{L/K}}
    &K_K^\times
  }
\end{displaymath}
in which the square is the case $E/L/L/K$ of the commutative square on
the left of \ref{eq:iota-delta-norms}, and the right hand triangle
commutes by the transitivity of norms. From this we see that
\begin{displaymath}
  N^K_{E/K}(\delta^E_{L/K}(\beta))=
  N^K_{L/K}(\delta^L_{L/K}(N^L_{E/L}(\beta)))
  =N^K_{L/K}(\delta^L_{L/K}(\eta_{E/L}(s))).
\end{displaymath}
Since $\eta_{E/L}\in L^\times$, lemma
\ref{lemma:random-norm-factoid} shows that
\begin{displaymath}
  N^K_{L/K}(\delta^L_{L/K}(\eta_{E/L}(s)))=N_{L/K}(\eta_{E/L}(s))
\end{displaymath}
and combining the previous equalities yields
\begin{equation}
  \label{eq:eta-subgroup}
  \eta_{E/K}(s)=N_{L/K}(\eta_{E/L}(s)) \mod N_{E/K}E^\times
\end{equation}
which is equivalent to
\begin{equation}
  \label{eq:norm-residue-subgroup}
  (a,E/L)=(N_{L/K}(a),E/K)
\end{equation}
for $a\in L^\times$.

The last compatibility asserts that the transfer
$\Ver:G_{E/K}^\ab\to G_{L/K}^\ab$ corresponds via the norm residue
symbol to the map
$K^\times/N_{E/K}E^\times\to L^\times/N_{E/L}E^\times$ induced by the
inclusion. From our point of view this is best understood in terms of
the ideas of the next section.

\section{Weil Groups}
\label{sec:weil-groups}

Usually the absolute Weil group of a local field is defined directly
in terms of the absolute Galois group; the formalism of the relative
Weil groups is an immediate consequence. In this section we show how
this formalism follows from our construction of the relative Weil
groups. This leads to a proof of Shafarevich's theorem.

\subsection{The formalism of Weil groups.}
\label{sec:Weil-formalism}

We fix a Galois extension $E/K$ with group $G_{E/K}$ and let $L$ be an
extension of $K$ in $E$ (not necessarily Galois over $K$). As before
$G_{E/L}$ is the Galois group of $E/L$, and similarly for $G_{L/K}$
when $L/K$ is Galois. Our aim in this section is to define morphisms
$i:W_{E/L}\to W_{E/K}$ and, when $L/K$ is Galois,
$p_{E/L}:W_{E/K}\to W_{L/K}$. Since $W_{L/K}$ by definition is only
defined up to inner automorphisms it will be necessary to ``rigidify''
it by a particular choice of ``realization'' $W(\alpha)$, which as
before is the set of pairs $(\beta,s)\in L_K^\times\times G_{L/K}$
satisfying \ref{eq:elements-of-W} with the composition law
\ref{eq:group-law-of-W}. To indicate the fields involved we write
$W_{L/K}(\alpha)$ for $W(\alpha)$, so a homomorphism such as
$p_{E/L}:W_{E/K}\to W_{L/K}$ should be understood as a homomorphism
$W_{E/K}(\alpha')\to W_{L/K}(\alpha)$ for particular choices (usually
implicit) of $\alpha$, $\alpha'$ appropriate to $L/K$ and $E/K$.

In fact everything we need is contained in the formulas
\ref{eq:eta-quotient} and \ref{eq:eta-subgroup}. The first of them
shows that $N^K_{E/L}:E_K^\times\to L_K^\times$ extends to a
homomorphism 
\begin{equation}
  \begin{split}
  p^K_{E/L}:W_{E/K}(\alpha)&\to W_{L/K}(N^K_{E/L}(\alpha))\\
  p^K_{E/L}(\beta,s)&=(N^K_{E/L}(\beta),\pi_{L/K}(s))
  \end{split}
  \label{eq:Weil-pi}
\end{equation}
making commutative the diagram
\begin{equation}
  \label{eq:Weil-galois-pi}
  \xymatrix{
    1\ar[r]&E^\times\ar[r]\ar[d]^{N_{E/L}}&W_{E/K}\ar[r]\ar[d]^{p^K_{E/L}}
    &G_{E/K}\ar[r]\ar[d]^{\pi_{E/L}}&1\\
    1\ar[r]&L^\times\ar[r]&W_{L/K}\ar[r]&G_{L/K}\ar[r]&1\\
  }.
\end{equation}
The transitivity of the norm shows that the maps $p$ just
constructed are transitive: 
\begin{equation}
  \label{eq:Weil-pi-transitive}
  p^K_{E/L}\circ p^K_{F/E}=p^K_{F/L}
\end{equation}
for a 3-fold extension $F/E/L/K$; here of course the groups $W_{F/K}$,
$W_{E/K}$ and $W_{L/K}$ are realized as $W_{F/K}(\alpha)$,
$W_{E/K}(N_{F/L}(\alpha)$ and $W_{L/K}(N_{F/L}(\alpha)$ for
appropriate $\alpha$; from now on we will not be explicit about this.

In the special case when $E/L/K$ is $L/K/K$ we identify
$W_{K/K}=K^\times$, and we will usually write $p_{L/K}$ for
$p^K_{L/K}:W_{L/K}\to K^\times$.  The definition of shows that the
diagram
\begin{equation}
  \label{eq:pi-and-eta}
  \xymatrix{
    W_{L/K}\ar[r]\ar[d]_{p_{L/K}}&G_{L/K}\ar[d]^{\eta_{L/K}}\\
    K^\times\ar[r]&K/N_{L/K}L^\times
  }
\end{equation}
commutes, where the horizontal maps are the natural projections. In
fact if $\alpha\in L_K^\times$ is used to define $W_{L/K}$,
$(\beta,s)\in W_{L/K}$ implies that
$\up{\sigma-1}\beta=\up{s-1}\alpha$, so that
\begin{displaymath}
  (p_{L/K}(\beta,s)\mod N_{L/K}L^\times)
  =(N_{L/K}(\beta)\mod N_{L/K}L^\times)=\eta_{L/K}(s).
\end{displaymath}
That \ref{eq:pi-and-eta} is commutative could be expressed by saying
that $p_{L/K}:W_{L/K}\to K^\times$ is a lifting of the inverse
reciprocity map. It could also be rephrased as saying that the
diagram
\begin{equation}
  \label{eq:pi-and-reciprocity}
  \xymatrix{
    W_{L/K}^\ab\ar[r]^{p_{L/K}^\ab}\ar[rd]&K^\times\ar[d]^{\theta_{L/K}}\\
    &G_{L/K}^\ab
  }
\end{equation}
is commutative.

Now \ref{eq:artin-symbol-subgroup} and the surrounding discussion
shows that $\delta^E_{L/K}:E_L^\times\to E_K^\times$ extends to a
homomorphism
\begin{equation}
  \label{eq:Weil-i}
  \begin{split}
    i^E_{L/K}:W_{E/L}(\alpha)&\to W_{E/K}(\gamma^E_{L/K}(\alpha))\\
    i^E_{L/K}(\beta,s)&=(\delta^E_{L/K}(\beta),s)
  \end{split}
\end{equation}
making commutative a diagram
\begin{equation}
  \label{eq:Weil-galois-i}
  \xymatrix{
    1\ar[r]&E^\times\ar[r]\ar@{=}[d]&W_{E/L}\ar[r]\ar[d]^{i^E_{L/K}}
    &G_{E/L}\ar[r]\ar[d]&1\\
    1\ar[r]&E^\times\ar[r]&W_{E/K}\ar[r]&G_{E/K}\ar[r]&1\\    
  }
\end{equation}
where the right vertical arrow is the canonical inclusion. For a
3-fold extension $F/E/L/K$ we have
\begin{equation}
  \label{eq:delta-transitive}
  \delta^F_{L/K}\circ\delta^F_{E/L}=\delta^F_{E/K}
\end{equation}
which shows that the $i$ are transitive:
\begin{equation}
  \label{eq:Weil-i-transitive}
  i^F_{L/K}\circ i^F_{E/L}=i^F_{E/K}.
\end{equation}
Finally the right hand diagram of \ref{eq:iota-delta-norms} implies
that for any 3-fold extension $F/E/L/K$ the diagram
\begin{equation}
  \label{eq:pi-i-commute-diagram}
  \xymatrix{
    W_{F/L}\ar[r]^{p^K_{F/E}}\ar[d]_{i^K_{L/K}}&W_{E/L}\ar[d]^{i^K_{L/K}}\\
    W_{F/K}\ar[r]_{p^K_{F/E}}&W_{E/K}
  }
\end{equation}
is commutative. This can be expressed by saying that if we use the
maps $i_{E/L}$ to identify $W_{E/L}$ with a subgroup of $W_{E/K}$,
then these identifications are compatible with the canonical
projection maps $p$ (assuming as always that compatible choices of
$\alpha$ are made in each place).

\subsection{The Shafarevich-Weil theorem.}
\label{sec:shafarevich-weil}

It is evident from the construction that $i^E_{L/K}:W_{E/L}\to
W_{E/K}$ is injective. It is also true that $p^K_{E/L}:W_{E/K}\to
W_{E/L}$ is surjective, although this is not completely obvious from
our construction.

\begin{prop}\label{prop:canonical-map-prop}
  The connecting homomorphism $\partial:G_{E/L}\to L^\times/NE^\times$
  arising from the diagram \ref{eq:Weil-galois-pi} is the inverse norm
  residue homomorphism $\eta_{E/L}$.
\end{prop}
\begin{demo} For $s\in G_{E/L}$, $\partial(s)$ is computed as follows:
  lift $s$ to an element $(\beta,s)\in W_{E/K}$; then
  $p_{E/L}(\beta,s)\in W_{L/K}$ lies in $L^\times\subset W_{L/K}$
  and $\partial(s)$ is the image of $p_{E/L}(\beta,s)$ in
  $L^\times/N_{E/L}E^\times$. The commutative diagram
  \ref{eq:pi-i-commute-diagram} applied to $E/L/L/K$ is
  \begin{displaymath}
      \xymatrix{
    W_{E/L}\ar[r]^{p_{E/L}}\ar[d]_{i_{L/K}}&W_{L/L}\ar[d]^{i_{L/K}}\\
    W_{E/K}\ar[r]_{p^K_{E/L}}&W_{L/K}
  }
  \end{displaymath}
  which shows that $\partial(s)$ can also be computed by lifting $s$
  to $(\beta,s)\in W_{E/L}$, and then $\partial(s)$ is the class of
  $p_{E/L}(\beta,s)\in W_{L/L}=L^\times$ in
  $L^\times/N_{E/L}E^\times$. Now $p_{E/L}(\beta,s)=N_{E/L}(\beta)$
  where $N_{E/L}$ is the norm $E_L^\times\to(L^\nr)^\times$, and if
  $\alpha\in E_L^\times$ is used to define $W_{E/L}$, $\beta$
  satisfies $\up{\sigma-1}\beta=\up{s-1}\alpha$, so
  $\partial(s)=\eta_{E/L}(s)$ by Dwork's formula.
\end{demo}

\begin{cor}\label{cor:weil-pi-surjective}
  For any successive extension $E/L/K$ of local fields with $E/K$ and
  $L/K$ Galois, the homomorphism $p^K_{E/L}:W_{E/K}\to W_{L/K}$ is
  surjective. 
\end{cor}
\begin{demo}
  Since $\partial=\eta_{E/L}$ is surjective, this follows from the
  snake lemma.
\end{demo}

\begin{lemma}\label{lemma:pi-and-transfer}
  For any Galois extension of local fields the transfer
  $\Ver:W_{L/K}^\ab\to L^\times$ is the composite
  \begin{displaymath}
    W_{L/K}^\ab\Xto{p_{L/K}^\ab}K^\times\to L^\times
  \end{displaymath}
  where we identify $W_{K/K}\simeq K^\times$, and the second map is
  the inclusion.
\end{lemma}
\begin{demo}
  Recall the definition of the transfer $G^\ab\to H^\ab$ for a
  subgroup $H\sset G$ of finite index: choose a section
  $\theta:H\setminus G\to G$ of the projection $G\to H\setminus G$,
  and for $s$, $t\in G$ define $x_{s,t}\in H$ by
  $\theta(Ht)s=x_{t,s}\theta(Hts)$; then
  \begin{displaymath}
    \Ver(s)=\prod_{t\in H\setminus G}x_{t,s}
  \end{displaymath}
  where the product is in $H^\ab$. In our case $G=W_{L/K}$ and
  $H=L^\times$ is abelian and normal in $G$. For $t\in G_{L/K}=G/H$ we
  take $\theta(t)=(\beta_t,t)$ for some appropriate
  $\beta_t\in L_K^\times$; then
  \begin{displaymath}
    \theta(t)(\beta,s)=(\beta_t,t)(\beta,s)=(\beta_t\cdot\up
    t\beta,ts)= 
    (\beta_t\up t\beta\beta_{ts}^{-1},1)\theta(ts)
  \end{displaymath}
  where the last equality makes sense since
  $\beta_t\up t\beta\beta_{ts}^{-1}\in L^\times$. For this $\theta$
  then we have $x_{t,s}=(\up s\beta,1)$, and then
  \begin{displaymath}
    \Ver(\beta,s)=\prod_{t\in G_{L/K}}(\beta_t\up
    t\beta\beta^{-1}_{ts},1)=(N_{L/K}(\beta),1)
  \end{displaymath}
  which proves the assertion.
\end{demo}

We can now prove the last compatibility for the norm residue symbol.

\begin{prop}\label{prop:transfer-compatibility}
  For any 2-fold extension $E/L/K$ of local fields with $E/K$ Galois,
  \begin{equation}
    \label{eq:transfer-compatibility}
    \xymatrix{
      K^\times\ar[r]\ar[d]_{\theta_{E/K}}&L^\times\ar[d]^{\theta_{E/L}}\\
      G_{E/K}^\ab\ar[r]^\Ver&G_{E/L}^\ab
    }
  \end{equation}
  is commutative, where the upper horizontal map is the inclusion.
\end{prop}
\begin{demo}
  By corollary \ref{cor:weil-pi-surjective}, the horizontal arrow in
  diagram \ref{eq:pi-and-reciprocity} is surjective. Applying this
  observation to the extensions $E/K$ and $E/L$, we see that it
  suffices to show that the diagrams
  \begin{equation}
    \label{eq:transfer-compatibility2}
    \xymatrix{
      W_{E/K}^\ab\ar[r]^{\Ver}\ar[d]&W_{E/L}^\ab\ar[d]\\
      G_{E/K}^\ab\ar[r]_{\Ver}&G_{E/L}^\ab
    }
    \qquad
    \xymatrix{
      W_{E/K}^\ab\ar[r]^{\Ver}\ar[d]_{p_{E/K}^\ab}
      &W_{E/L}^\ab\ar[d]^{p_{E/L}^\ab}\\
      K^\times\ar[r]_{\mathit{incl}}&L^\times
    }
  \end{equation}
  are commutative. The first just expresses the functoriality of the
  transfer. The second can be embedded in the diagram
  \begin{displaymath}
    \xymatrix{
      W_{E/K}^\ab\ar[r]^{\Ver}\ar[d]_{p_{E/K}^\ab}
      &W_{E/L}^\ab\ar[d]^{p_{E/L}^\ab}\ar@/^.7cm/[rdd]^\Ver\\
      K^\times\ar[r]\ar@/_.5cm/[rrd]_i&L^\times\ar[rd]_j\\
      &&E^\times      
    }
  \end{displaymath}
  in which $i$ and $j$ are the inclusions. By lemma
  \ref{lemma:pi-and-transfer} $i\circ p^\ab_{E/K}$ and
  $j\circ p^\ab_{E/L}$ are the transfer homomorphisms for the
  subgroups $E^\times\subset W_{E/K}$ and $E^\times\subset W_{E/L}$
  respectively. Since the transfer is transitive for chains of
  subgroups, the outside square in this diagram is commutative, and
  since $L^\times\to E^\times$ is injective, the inside square is as
  well.
\end{demo}

In terms of the norm residue symbol, the commutativity of
\ref{eq:transfer-compatibility} says that
\begin{equation}
  \label{eq:transfer-compatibility3}
  \Ver(a,E/K)=(a,E/L)
\end{equation}
for $a\in K^\times$. Proposition \ref{prop:transfer-compatibility} is
less elementary than the previous compatibilities since the
homomorphism $G^\ab_{E/K}\to G^\ab_{E/L}$ is not induced by a
homomorphism $G_{E/K}\to G_{E/L}$.

We can now prove Shafarevich's theorem
\cite{shafarevich:1946}. Suppose that $E/L/K$ are extensions of local
fields with $E/K$ and $L/K$ Galois. We define a map
\begin{equation}
  \label{eq:shafarevich-map}
  \Sh:G_{E/K}\to W_{L/K}/N_{E/L}E^\times
\end{equation}
as follows: for $s\in G_{E/K}$ choose
$x\in W_{E/K}$ mapping to $s$ under the natural projection $W_{E/K}\to
G_{E/K}$. Then $p^E_{L/K}(x)\in W_{L/K}$ is well-defined up to a
factor in $N_{E/L}E^\times\subset L^\times\subset W_{L/K}$, and we
define $\Sh(s)$ to be the class of $p^E_{L/K}(x)$ in
$W_{L/K}/N_{E/L}E^\times$. 

\begin{thm}[Shafarevich]\label{thm:shafarevich}
  Suppose that $E/L/K$ are extensions of local fields with $E/K$ and
  $L/K$ Galois. The diagram
  \begin{equation}
    \label{eq:shafarevich}
    \xymatrix{
      1\ar[r]&G_{E/L}\ar[r]\ar[d]^{\eta_{E/L}}&G_{E/K}\ar[r]\ar[d]^\Sh
      &G_{L/K}\ar[r]\ar@{=}[d]&1\\
      1\ar[r]&L^\times/N_{E/L}E^\times\ar[r]&W_{L/K}/N_{E/L}E^\times\ar[r]
      &G_{L/K}\ar[r]&1
    }
  \end{equation}
  is commutative. In particular if $E/L$ is abelian, $\Sh$ is an
  isomorphism.
\end{thm}
\begin{demo}
  The commutativity of the right hand square is clear from the
  construction. For the one on the left, we identify
  $W_{E/L}=W_{E/L}(\alpha)$ for some $\alpha\in E_L^\times$, and
  $W_{E/K}=W_{E/K}(\gamma^E_{L/K}(\alpha))$ as in our construction of
  $i^E_{L/K}:W_{E/L}\to W_{E/K}$. Any lifting of an $s\in G_{E/L}$ to
  $W_{E/K}$ lies in the image of $i^E_{L/K}$, so we may suppose the
  lifting is $(\delta^E_{L/K}(\beta),s)$, with $\beta\in L_K^\times$
  satisfying $\up{\sigma-1}\beta=\up{s-1}\alpha$. Now the diagram
  \ref{eq:pi-i-commute-diagram} in the case when $F/E/L/K$ is
  $E/L/L/K$ is
  \begin{displaymath}
      \xymatrix{
    W_{E/L}\ar[r]^{p^L_{E/L}}\ar[d]_{i^E_{L/K}}&W_{L/L}\ar[d]^{i^L_{L/K}}\\
    W_{E/K}\ar[r]_{p^K_{E/L}}&W_{L/K}
  }
  \end{displaymath}
  which shows that $\Sh(s)$ can be identified with the image of
  \begin{displaymath}
    p^L_{E/L}(\beta,s)=(N_{E/L}(\beta),1)\in W_{L/K}.
  \end{displaymath}
  Since $N_{E/L}(\beta)=\eta_{E/L}(s)$ the left hand square is
  commutative. Finally if $E/L$ is abelian $\eta_{E/L}$ is an
  isomorphism, hence so is $\Sh$.
\end{demo}


\begin{thebibliography}{99}
\bibitem{bourbaki:csa} N. Bourbaki, \emph{Alg\`ebre, Ch. 8. Modules
    et anneaux semi-simples}, Springer 2012.
\bibitem{bourbaki:hom} N. Bourbaki, \emph{Alg\`ebre,
    Ch. 11. Alg\'ebre homologique} , Springer 2007.
\bibitem{deligne:1972} P. Deligne, \emph{Les constantes des
    \'equations fonctionelles des fonctions $L$}, in \emph{Modular
    Functions of One Variable II}, Springer LNM \textbf{349} (1973)
  pp.~501--597.
\bibitem{demazure:1972} M. Demazure, \emph{Lectures on $p$-Divisible
    Groups}, Springer LNM \textbf{302}, 1972.
\bibitem{dwork:1958}B. Dwork, \emph{Norm residue symbol in local
    number fields}, Hamb. Abh. \textbf{22} (1958) pp.~180--190.
\bibitem{katz:1979} N. Katz, \emph{Slope Filtration of
    $F$-crystals.}, in \emph{Journ\'ees de G\'eom\'etrie Algebrique
    de Rennes I}, Ast\'erisque \textbf{63} (1979), SMF.
\bibitem{manin:1963} Yu. Manin, \emph{The Theory of Commutative Formal
    Groups over Fields of Finite Characteristic}, Russ. Math. Surveys
  \textbf{18} pp.~1--83.
\bibitem{morava:1979} J. Morava, \emph{The Weil group as an
    automorphism group}, in \emph{Journ\'ees de G\'eom\'etrie
    Algebrique de Rennes I}, Ast\'erisque \textbf{63} (1979), SMF.
\bibitem{serre:1967} J.-P. Serre, \emph{Local Class Field Theory},
  in \emph{Algebraic Number Theory}, J. Cassels and A. Fr\"ohlich
  (eds), Academic Press (1967) pp.~128--161.
\bibitem{serre:1968} J.-P. Serre, \emph{Corps locaux}, Hermann
  1968.
\bibitem{shafarevich:1946} I. R. Shafarevich, \emph{On Galois groups
    of $p$-adic fields}, Doklady Acad. Sc. USSR (N.S.) \textbf{53}
  (1946) pp.~4--5.
\bibitem{tate:arcata} J. Tate, \emph{Number Theoretic Background},
  in \emph{Automorphic Forms, Representations and $L$-functions},
  PSPM \textbf{33} (1979), part 2 pp.~3--26.
\end{thebibliography}
\end{document}